\documentclass[12pt,a4paper]{amsart}
\usepackage{amssymb,euler}
\def\COM{\mathbb C}
\def\REA{\mathbb R}

\def\J{J}
\def\K{K}
\def\O{O}
\def\T{L}
\def\N{k}
\def\R{\tau}
\def\A{\Delta}
\def\C{C}
\def\h{h}
\def\m{m}
\def\p{p}
\def\q{q}
\def\I{\mathsf i}

\newtheorem*{theorem}{Theorem}
\newtheorem{lemma}{Lemma}
\newtheorem*{corollary}{Corollary}
\theoremstyle{remark}
\newtheorem*{ack}{Acknowledgments}

\title{Proof of the volume conjecture for torus knots}
\author{R.~M.~Kashaev}
\address{Steklov Mathematical Institute at St. Petersburg,
Fontanka 27, St. Petersburg 191011, Russia}
\curraddr{Helsinki Institute of Physics, P.O. Box 9, FIN-00014
University of Helsinki, Finland}
\email{kashaev@pdmi.ras.ru}
\author{O.~Tirkkonen}
\address{Nokia Research Center, P.O. Box 407, FIN-00045, Nokia Group, Finland }
\email{olav.tirkkonen@nokia.com}
\date{December 1999}

\begin{document}
\begin{abstract}
The volume conjecture, formulated recently by H. Murakami and J. Murakami,
is proved for the case of torus knots. 
\end{abstract}
\maketitle
\section{Introduction}

In the recent paper \cite{Murakami} H. Murakami and J. Murakami
showed that the ``quantum dilogarithm'' knot invariant, introduced
in \cite{Kas1,Kas2}, is a special case of the
colored Jones invariant (polynomial)
associated with the quantum group $SU(2)_\q$. Using the 
connection of the quantum dilogarithm invariant with the hyperbolic volume of
knot's complement, conjectured in 
\cite{Kas3}, they also proposed the
``volume conjecture'': for any knot the above mentioned specification
of the colored Jones invariant in a certain limit gives the simplicial volume
(or Gromov norm) of the knot. It is remarkable that this conjecture
implies that a knot is trivial if and only if its colored Jones invariants are 
trivial, see \cite{Murakami}.

The purpose of this paper is to prove the volume conjecture for the case
of torus knots. To formulate our result, let us first
recall the form of the colored 
Jones invariant for torus knots \cite{Morton,RossoJones}.

The colored Jones 
invariant $\J_{\K,\N}(\h)$ of a framed knot $\K$ is a Laurent polynomial in 
$\q=e^\h$ depending on the `color' $\N$, the dimension of
a $SU(2)_\q$-module. Let $\m,\p$ be mutually prime positive integers.
Denote $\T\equiv \O_{\m,\p}$ the $(\m,\p)$ torus knot obtained 
as the $(\m,\p)$ cable about the unknot with zero
framing (see \cite{Morton} for the precise definition).
Then the colored Jones invariant of $\T$
has the following explicit form:
\begin{equation}\label{torus1}
2 \sinh(\N \h/2)\frac{\J_{\T,\N}(\h)}{\J_{\O,\N}(\h)}=
\sum_{\epsilon=\pm1}\sum_{r=-(\N-1)/2}^{(\N-1)/2}
\epsilon e^{\h\m\p r^2+\h r(\m+\epsilon \p)+\epsilon \h/2},
\end{equation}
where $\O$ is the unknot with zero framing, and 
\[
\J_{\O,\N}(\h)=\sinh(\N \h/2)/\sinh(\h/2).
\]
According
to the H. Murakami and J. Murakami's result \cite{Murakami}, 
the quantum dilogarithm invariant $\langle \T\rangle_{\N}$
is the following specification 
of the colored Jones invariant (up to a multiple of an $\N$-th root of unity):
\begin{equation}\label{hyperbolic}
\langle \T\rangle_{\N}\equiv 
\lim_{\h\to 2 \pi \I/\N}\frac{\J_{\T,\N}(\h)}{\J_{\O,\N}(\h)}.
\end{equation}
In what follows, we shall call this as ``hyperbolic specification''.
Our result describes the asymtotic 
expansion of $\langle \T\rangle_{\N}$
when $\N\to\infty$.

\begin{theorem}
The hyperbolic specification (\ref{hyperbolic})
of the colored Jones invariant for the $(\m,\p)$ torus knot $\T$
has the following asymtotic expansion at large $\N$:
\begin{equation}\label{mainrez}
\langle \T\rangle_{\N}
e^{\frac{\I\pi}{2\N}\left(\frac\m\p+\frac\p\m\right)}
=\sum_{j=1}^{\m\p-1}\langle \T\rangle_{\N}^{(j)}+
\langle \T\rangle_{\N}^{(\infty)},
\end{equation}
where
\begin{equation}\label{mainrez1}
\langle \T\rangle_{\N}^{(j)}
=2(2\m\p/\N)^{-3/2}e^{\I\frac\pi 4}
(-1)^{(\N-1)j}e^{-\I\frac{\pi \N j^2}{2\m\p}}
j^2 \sin(\pi j/m)\sin(\pi j/p),
\end{equation}
and
\begin{equation}\label{mainrez2}
\langle \T\rangle_{\N}^{(\infty)}
=\frac14e^{\I\pi \N \m \p/2}
\sum_{n\ge 1}\frac1{n!}\left(\frac{\I\pi}{2 \N\m\p}\right)^{n-1}
\left.\frac{\partial^{2n}(x \R_\T(x))}{\partial x^{2n}}\right\vert_{x=0},
\end{equation}
see formula~(\ref{rtorsion}) below for the definition of the function
$\R_\T(x)$.
\end{theorem}
From eqns~(\ref{mainrez})--(\ref{mainrez2}) it is easily seen that
$|\langle \T\rangle_{\N}|\sim\N^{3/2}$, $\N\to\infty$.
\begin{corollary}
The volume conjecture holds true for all torus knots, i.e.
\[
\lim_{\N\to \infty}\N^{-1}\log|\langle \T\rangle_{\N}|=0.
\]
\end{corollary}
In the next section we prove the Theorem by using an integral representation
for the Gaussian
sum in formula~(\ref{torus1}).

\begin{ack}
We are grateful to D. Borisov, T. K\"arki and H. Murakami for discussions. 
R.K. thanks L.D. Faddeev for his constant support in this work.
The work of R.K. is supported by Finnish Academy, and in part 
by RFFI grant 99-01-00101. 
\end{ack}
\section{Proof of the Theorem}
 To begin with, define the following function:
\begin{equation}\label{rtorsion}
\R_\T(z)\equiv 2\sinh(mz)\sinh(pz)/\sinh(mpz).
\end{equation}
It is related to the Alexander polynomial of the knot $\T$,
\[
\A_\T(t)\equiv\frac{(t^{\m\p/2}-t^{-\m\p/2})(t^{1/2}-t^{-1/2})}
{(t^{\m/2}-t^{-\m/2})(t^{\p/2}-t^{-\p/2})},
\]
through the formula
\[
\R_\T(z)=2\sinh(z)/
\A_\T(e^{2z}).
\]
According to the result of Milnor \cite{Milnor} and Turaev \cite{Turaev},
the function $\R_\T(z)$ describes the Reidemeister torsion of 
the knot complement.
\begin{lemma}
For any real $\phi$, satisfying the condition 
$\Re\h e^{-2\I\phi}>0$, formula (\ref{torus1}) 
has the following integral representation
\begin{equation}\label{integral1}
2\sinh(\N \h/2)
\frac{\J_{\T,\N}(\h)}{\J_{\O,\N}(\h)}=\sqrt{\frac{\m\p}{\pi\h}}
e^{-\frac\h 4\left(\frac\m\p+\frac\p\m\right)}
\int_{\C_\phi}dz\, e^{\m\p\left(\N z-\frac{z^2}\h\right)}\R_\T(z),
\end{equation}
where the integration path $\C_\phi$ is the image of the real line 
under the mapping
\begin{equation}\label{contour}
\REA\ni x\mapsto xe^{\I\phi}\in\C_\phi\subset\COM,
\end{equation}
with the induced orientation.
\end{lemma}
\begin{proof}
First note that for any complex $\h\ne 0$ and any complex $w$, the 
following Gaussian integral formula holds:
\begin{equation}\label{gauss}
\sqrt{\pi\h}e^{\h w^2}=\int_{\C_\phi}dz\, e^{-z^2/\h+2wz},
\end{equation}
where the choice of the integration path $\C_\phi$, described in the 
formulation of the theorem, is dictated by the 
convergence condition of the integral, and the square root is the analytical 
continuation from positive values of $\h$.
Now, starting from the right hand side of eqn~(\ref{torus1}), 
collect the terms, containing the summation variable $r$, 
into a complete square:
\[
e^{\frac\h 4\left(\frac\m\p+\frac\p\m\right)}\,\mathrm{r.h.s.}(\ref{torus1})=
\sum_{\epsilon=\pm1}\epsilon\sum_{r=0}^{\N-1}
e^{\h\m\p \left(r-\frac{\N-1}2+\frac{\m+\epsilon\p}{2\m\p}\right)^2}
\]
--- now formula (\ref{gauss}) can be applied  to the $r$-dependent 
exponential~---
\[
=\sum_{\epsilon=\pm1}\epsilon\sum_{r=0}^{\N-1}\frac1{\sqrt{\pi\h\m\p}}
\int_{\C_\phi}dz\, e^{-\frac{z^2}{\h\m\p}+
z(2r-\N+1+\p^{-1}+\epsilon\m^{-1})}
\]
--- with subsequent evaluation of the summations~---
\[
=\frac2{\sqrt{\pi\h\m\p}}\int_{\C_\phi}dz\, e^{-\frac{z^2}{\h\m\p}+z/\p}
\sinh(\N z)\sinh(z/m)/\sinh(z)
\]
--- the exponential $\exp(z/\p)$ in the integrand, being multiplied by an
odd function of $z$ (w.r.t. $z\leftrightarrow -z$), 
can by replaced by it's odd part ---
\[
=\frac2{\sqrt{\pi\h\m\p}}\int_{\C_\phi}dz\, e^{-\frac{z^2}{\h\m\p}}
\sinh(\N z)\sinh(z/m)\sinh(z/p)/\sinh(z)
\]
--- reversing the previous argument, replace $\sinh(\N z)$ by
an exponential~---
\[
=\frac2{\sqrt{\pi\h\m\p}}\int_{\C_\phi}dz\, e^{-\frac{z^2}{\h\m\p}+\N z}
\sinh(z/\m)\sinh(z/p)/\sinh(z)
\]
--- and rescale the integration variable ($z\to z\m\p$)~---
\[
=\sqrt{\frac{\m\p}{\pi\h}}
\int_{\C_\phi}dz\, e^{\m\p\left(\N z-\frac{z^2}\h\right)}\R_\T(z) 
\]
with notation (\ref{rtorsion}) being used.
\end{proof}
 Representation (\ref{integral1})
is similar to Rozansky's formula (2.2) from \cite{Roz}, though the latter is 
only a shorthand for the power series expansion.
\begin{lemma}
The hyperbolic specification (\ref{hyperbolic}) of the colored Jones 
invariant for torus knots has the integral representation
\begin{equation}\label{hypint1}
2\langle \T\rangle_{\N}=(\m\p\N/2)^{3/2}
e^{-\frac{\I\pi}{2\N}\left(\frac\m\p+\frac\p\m+\frac\N2\right)}
\int_{\C_\phi}dz\,  e^{\pi\m\p\N(z+\frac\I2 z^2)}z^2\R_\T(\pi z),
\end{equation}
where integration path $\C_\phi$ is defined in (\ref{contour}) with
$0<\phi<\pi/2$.
\end{lemma}
\begin{proof} The left hand side of eqn~(\ref{integral1}) vanishes
at $\h=2\pi\I/\N$ due to the factor $\sinh(\N \h/2)$. This means that
the integral in the right hand side vanishes as well. So, differentiating
simultaneously the $\sinh$-function in the left hand side and the 
integral in the right of eqn~(\ref{integral1}) with respect to $\h$, then
putting $\h=2\pi\I/\N$, and rescaling the integration variable by $\pi$,
we rewrite the result in the form
of eqn~(\ref{hypint1}).
\end{proof}

\begin{proof}[Proof of the Theorem]
At large $\N$ one can use the steepest descent method for evaluation
the integral in (\ref{hypint1}). The only stationary point
at $z=\I$ is separated from the integration path by a finite number of 
poles of the function
$\R_\T(\pi z)$ which are located at $z_j\equiv\I j/\m\p$, $0<j<\m\p$. 
Thus, 
taking into account convergence at infinity, we can shift path $\C_\phi$
by imaginary unit and add integration along a closed contour
encircling points $z_j$ in the counterclockwise direction. The integration
along the shifted path $\C_\phi$ can be transformed by the change of the 
integration variable $z\to z+\I$: 
\begin{multline*}
\int_{\I+\C_\phi}dz\,  e^{\pi\m\p\N(z+\frac\I2 z^2)}z^2\R_\T(\pi z)\\
=
 e^{\I\pi\m\p\N/2}\int_{\C_\phi}dz\,  e^{\I\pi\m\p\N z^2/2}(z+\I)^2
\R_\T(\pi z+\I\pi)\\
= -2\I e^{\I\pi\m\p\N/2}\int_{\C_\phi}dz\,  
e^{\I\pi\m\p\N z^2/2}z
\R_\T(\pi z),
\end{multline*}
where in the last line we have used the (quasi) periodicity property of the 
$\sinh$-function, the fact that $\m$ and $\p$ are mutually prime, and
disregarded the odd terms with respect to the sign change 
$z\leftrightarrow-z$. Now, the obtained formula straightforwardly leads
to the asymptotic power series $\langle \T\rangle_{\N}^{(\infty)}$ in 
eqn~(\ref{mainrez}) through the Taylor series expansion of the 
function $z\R_\T(\pi z)$ at $z=0$, and evaluation of the Gaussian integrals. 
The other terms in eqn~(\ref{mainrez}) come from the 
evaluation of the contour integral by the residue method.
\end{proof}

\end{document}